\documentclass[10pt]{amsart}
\calclayout
\usepackage[T1]{fontenc}
\usepackage{lmodern}
\usepackage{microtype}
\usepackage{amsmath,amssymb,amsthm,mathtools}
\IfFileExists{stmaryrd.sty}{\usepackage{stmaryrd}}{}
\usepackage{aliascnt}
\usepackage{enumitem}
\usepackage{array,booktabs,longtable,multirow}
\usepackage{xcolor}
\usepackage[colorlinks=true,linkcolor=blue!45!black,citecolor=blue!45!black,urlcolor=blue!45!black]{hyperref}
\usepackage[nameinlink,noabbrev]{cleveref}

\newtheorem{theorem}{Theorem}[section]
\newtheorem{claim}{Claim}
\newaliascnt{proposition}{theorem}
\newtheorem{proposition}[proposition]{Proposition}
\aliascntresetthe{proposition}
\newaliascnt{lemma}{theorem}
\newtheorem{lemma}[lemma]{Lemma}
\aliascntresetthe{lemma}
\newaliascnt{question}{theorem}

\aliascntresetthe{question}
\newaliascnt{corollary}{theorem}

\aliascntresetthe{corollary}
\theoremstyle{definition}
\newaliascnt{definition}{theorem}
\newtheorem{definition}[definition]{Definition}
\aliascntresetthe{definition}
\newaliascnt{construction}{theorem}

\aliascntresetthe{construction}
\theoremstyle{remark}
\newaliascnt{remark}{theorem}
\newtheorem{remark}[remark]{Remark}
\aliascntresetthe{remark}

\newcommand{\N}{\omega}
\newcommand{\Tot}{\operatorname{Tot}}
\newcommand{\BLR}{\operatorname{BLR}}

\begin{document}

\title{Some Variants of Function Domination}

\author{Xuanheng Zhao}
\address{School of Mathematics, Nanjing University,
Nanjing, Jiangsu 210093, People's Republic of China}
\email{xuanheng21@gmail.com}


\begin{abstract}
We investigate whether classical recursion-theoretic notions defined by
function domination remain unchanged when success is
required only on sets of prescribed asymptotic density. We establish robustness
for many such notions and uncover unexpected connections
between familiar degree properties by comparing upper and lower density. In
contrast, pdomination, as defined by Chong, Hoi, Stephan and Turetsky, is not robust under density weakenings.
We also answer a question of Kjos-Hanssen and Nies by characterizing
superhighness as domination with a threshold boundedly limit recursive in the
oracle. As an application, we introduce supermaximal sets, a strengthening of
maximality whose degrees are exactly the superhigh
r.e. degrees.
\end{abstract}

\maketitle
\setcounter{tocdepth}{2}

\section{Introduction and preliminaries}\label{sec:1}

Function domination is one of the principal ways in which growth rates encode
the computational strength of a Turing degree.  For total functions
$f,g\colon\N\to\N$, we say that $f$ \emph{dominates} $g$, and write
$g\ll f$, if $(\forall^\infty n)\;g(n)\leq f(n)$.
A function is \emph{recursively dominated} if some recursive function
dominates it.  Following Miller and Martin~\cite{MR228341}, a degree is
\emph{hyperimmune-free} if every function in it is recursively dominated,
and it is \emph{hyperimmune} otherwise.  In the opposite direction, Martin's
domination theorem~\cite{MR224469} identifies the high degrees:
\[
  A'\geq_T\emptyset''
  \quad\Longleftrightarrow\quad
  (\exists f\leq_TA)(\forall g\text{ recursive})\;g\ll f.
\]
These two classical notions already suggest that the placement of the
quantifiers, rather than growth rate alone, carries the degree-theoretic
content.

This paper asks how much of that content survives when eventual or everywhere
success is weakened to success on a large set of inputs.  The test is
asymptotic density.  For $C\subseteq\N$ and $n>0$, let $\rho_n(C)=\frac{|C\cap[0,n)|}{n}$,
and define its lower and upper asymptotic densities by $\underline\rho(C)=\liminf_{n\to\infty}\rho_n(C)$ and $\overline\rho(C)=\limsup_{n\to\infty}\rho_n(C)$.
If these two quantities are equal, their common value is denoted by
$\rho(C)$.

\begin{definition}\label{def:coarse-domination}
Given functions $f,g\colon\N\to\N$, we say that $f$ \emph{coarsely
dominates} $g$, written $g\ll_c f$, if
$\rho(\{m:f(m)\geq g(m)\})=1$.
A function $f$ is \emph{coarsely recursively bounded} if it is coarsely
dominated by some recursive function.
\end{definition}

Every recursively dominated function is coarsely recursively bounded, but
the converse fails, which follows from a simple diagonal argument.

Our main concern is \emph{robustness}.  A definition
based on domination, agreement, or tracing is robust under a density
weakening when the weakened condition defines the same degree class.
For many of the notions considered here, even positive lower density
retains the content of eventual or everywhere success.  In investigating
this robustness, we also uncover some unexpected connections: replacing
upper density by lower density in suitable formulations relates familiar
notions whose usual definitions look quite different.

\Cref{sec:density-degrees} begins by showing that coarse hyperimmunity is
exactly ordinary hyperimmunity, and that
highness is preserved when eventual domination is weakened to domination
on a set of positive lower density.  These
results also place hyperimmunity and highness into a common scheme:
\[
\begin{aligned}
\mathbf a\text{ is hyperimmune}
&\iff (\exists f\leq_T\mathbf a)(\forall\text{ recursive }g)\;
  \overline\rho\bigl(\{n:f(n) \ge g(n)\}\bigr)>0,\\
\mathbf a\text{ is high}
&\iff (\exists f\leq_T\mathbf a)(\forall\text{ recursive }g)\;
  \underline\rho\bigl(\{n:f(n)\geq g(n)\}\bigr)>0.
\end{aligned}
\]
In the upper-density characterization, positive upper density can be
strengthened to upper density one.  We also obtain an upper-density-one
version of the infinitely-often-equal characterization of hyperimmunity
(see Monin and Nies~\cite[Proposition~II.17]{MR3538473}).

Array recursiveness is also preserved when eventual domination is replaced
by domination on a set of positive lower density.  The corresponding upper-density
condition instead characterizes $\emptyset'\nleq_TA$.
For $A\leq_T\emptyset'$, lowness admits an analogous lower-density
characterization using a single $f\leq_T\emptyset'$ to bound all partial
$A$-recursive functions, with undefined inputs counted as correct.  Its upper-density version characterizes
incompleteness.

Recursive semi-traceability is unchanged if agreement at infinitely many
inputs is replaced by agreement on a set of upper density one.
Positive upper density yields the same class; positive lower density
instead forces $A$ to be recursive.

Strong jump traceability is likewise preserved when, for every recursive
order bound and every
partial $A$-recursive function, an r.e. trace obeying the prescribed bound
is required to be correct only on a set of positive lower density, again
counting undefined inputs as correct.  The upper-density analogue
requires a different treatment of partiality: we spread each input over a
long finite block to account for sparse domains, and count only genuine
trace hits at defined inputs.  Requiring that, for every order bound and
every partial $A$-recursive function $\psi$ with infinite domain, there be an
r.e. trace obeying that bound whose genuine hits for $\psi^{\mathrm{blk}}$
have positive upper density characterizes the non-DNR degrees.

Almost-everywhere domination provides a measure-theoretic instance of
robustness.  Using its equivalence with uniformly almost-everywhere and
positive-measure domination due to Kjos-Hanssen, Miller, and
Solomon~\cite{KMS12}, we show that positive lower density suffices in both
the uniform measure-one and the positive-measure formulations.  The uniform formulation requires one
$f\leq_TA$ to dominate every total $g\leq_TX$ on a set of positive lower
density, for almost every $X\in2^\N$.  Replacing lower density by upper
density here yields precisely the hyperimmune degrees.

\Cref{sec:bounded-threshold} characterizes superhighness through domination
with effective thresholds. As an application to r.e. degrees, we study
an effective strengthening of maximality for r.e. sets.
Mohrherr~\cite{Mohrherr86} defined a set $A$ to be
\emph{superhigh} when
$\emptyset''\leq_{\mathrm{tt}}A'$.
Superhighness occurs naturally in work on almost-everywhere domination and
algorithmic randomness~\cite{Simpson07,KHN09}. Kjos-Hanssen and
Nies~\cite{KHN09} asked for a domination characterization of superhighness.

\begin{definition}[Cole and Simpson~{\cite{ColeSimpson07}}]
\label{def:blr}
A total function $h\colon\N\to\N$ is \emph{boundedly limit recursive in
$A$} if there are an $A$-recursive function $H(n,s)$ and a recursive
function $r$ such that $ h(n)=\lim_sH(n,s)$ and 
$|\{s:H(n,s)\neq H(n,s+1)\}|<r(n)$
for every $n$.  We write
\[
  \BLR(A)=\{h:h\text{ is boundedly limit recursive in }A\}.
\]
\end{definition}

We answer their question in \Cref{thm:superhigh-threshold} (with Liang Yu): a set $A$ is
superhigh if and only if there are $f\leq_TA$ and $h\in\BLR(A)$ such that,
for the standard enumeration $(\Phi_e)_{e\in\N}$ of the partial recursive
functions, whenever $\Phi_e$ is total,
$(\forall n\geq h(e))\;\Phi_e(n)\leq f(n)$.
Thus superhighness is characterized by domination together with
bounded-limit control of when domination begins.

An application is that the degrees containing \emph{supermaximal}
sets are exactly the superhigh r.e. degrees
(\Cref{thm:superhigh-supermaximal}). This gives an analogue of Martin's
characterization of the high r.e. degrees as those containing maximal
sets~\cite{MR224469}. Recall that an r.e. coinfinite set $A$ is maximal
if its complement is cohesive: for each $e$, all sufficiently large
elements of $\overline A$ agree on membership in $W_e$. We call $A$
\emph{supermaximal} if these homogeneity thresholds can be chosen as the
values of a single function $h\in\BLR(A)$.
Since there are high r.e.
degrees that are not superhigh~\cite{Mohrherr86}, this defines a proper
subclass of the maximal sets.

\Cref{thm:superhigh-density} returns to density robustness. Superhighness
is also equivalent to the existence of $f\leq_TA$ for which each total
$\Phi_e$ is dominated on a set of density one, with convergence controlled
by a two-variable function $m\in\BLR(A)$: for every $k\geq1$ and
$N\geq\max\{1,m(e,k)\}$,
$\rho_N(\{n:\Phi_e(n)\leq f(n)\})\geq1-2^{-k}$.
For any fixed rational $0<\alpha<1$, it is already sufficient to require
these proportions to be at least $\alpha$ beyond a threshold in $\BLR(A)$.

\Cref{sec:ldpd-collapse} gives a contrasting example using pdomination,
introduced by Chong, Hoi, Stephan, and
Turetsky~\cite{ChongHoiStephanTuretsky15}. We consider two weakenings:
counting inputs outside the test domain as correct, and measuring successful
domination relative to that domain. Both collapse to the class of all sets
even at density one. A single
partial recursive function witnesses both density-one conditions. These
results show that density robustness cannot be taken for granted.

\subsection{Preliminaries}

We use standard notation from recursion theory and identify a Turing degree
with any representative when no confusion can arise.
For background, see Soare~\cite{Soare87}, Nies~\cite{Nies09}, and Downey and
Hirschfeldt~\cite{DowneyHirschfeldt10}.
Unless explicitly
declared partial, all functions have domain~$\N$. For sets $A$ and $B$, the
notations $A\leq_{\mathrm{tt}}B$ and $A\leq_{\mathrm{wtt}}B$ denote
truth-table and weak truth-table reducibility, respectively. 

We begin with a rapidly growing block system used in the density conversions
in \Cref{sec:density-degrees}. Its essential feature is that each block occupies
an asymptotically full proportion of the initial segment ending at its right
endpoint.

\begin{definition}
Let
$(q_n)_{n\in\N}$ and $(\mathcal Q_n)_{n\in\N}$ be given by $q_0=2$, $q_{n+1}=q_n^2$, and
\begin{equation}\label{eq:quadratic-blocks}
  \mathcal Q_n=[q_n,q_{n+1}).
\end{equation}
\end{definition}

\begin{lemma}
\label{lem:long-block-density}
For every $C\subseteq\N$,
\[
  \left|\rho_{q_{n+1}}(C)
    -\frac{|C\cap\mathcal Q_n|}{|\mathcal Q_n|}\right|
  \leq\frac1{q_n}\longrightarrow0.
\]
Consequently:
\begin{enumerate}[label=\textnormal{(\roman*)}]
\item If $\underline\rho(C)=\delta>0$ and $0\leq\eta<\delta$, then
$|C\cap\mathcal Q_n|>\eta|\mathcal Q_n|$ for all sufficiently large~$n$.
\item If $\mathcal Q_n\subseteq C$ for infinitely many~$n$, then
$\overline\rho(C)=1$.
\item If $C\cap\mathcal Q_n=\varnothing$ for infinitely many~$n$, then
$\underline\rho(C)=0$.
\end{enumerate}
\end{lemma}

\begin{proof}
Since $[0,q_{n+1})=[0,q_n)\cup\mathcal Q_n$ is a disjoint union,
\[
\begin{aligned}
  \left|\rho_{q_{n+1}}(C)
    -\frac{|C\cap\mathcal Q_n|}{|\mathcal Q_n|}\right|
  =\frac{q_n}{q_{n+1}}
    \left|\rho_{q_n}(C)
      -\frac{|C\cap\mathcal Q_n|}{|\mathcal Q_n|}\right|
  \leq\frac{q_n}{q_{n+1}}=\frac1{q_n}.
\end{aligned}
\]
The inequality holds because both proportions inside the absolute value
belong to $[0,1]$. Since $q_n\to\infty$, the three consequences follow
by taking lower or upper limits.
\end{proof}

\begin{definition}
For any partial function $\psi$, define its \emph{block spreading}
$\psi^{\mathrm{blk}}$ by
\begin{equation*}
  \psi^{\mathrm{blk}}(m)\simeq
  \begin{cases}
    \psi(0),&m<q_0,\\
    \psi(n),&m\in\mathcal Q_n.
  \end{cases}
\end{equation*}
\end{definition}

For total $\psi$, block spreading preserves weak truth-table
reducibility to a set $A$.
Applied to a positive, nondecreasing, unbounded recursive function, block
spreading preserves these properties as well. Together with
\Cref{lem:long-block-density}, it is the common mechanism behind the proofs in \Cref{sec:density-degrees}.

We next fix the effective indexing used for superhighness.  Let
$(\Phi_e)_{e\in\N}$ be a standard enumeration of the partial recursive
functions, and put
$\Tot=\{e:\Phi_e\text{ is total}\}$,
$\mathrm{INF}=\{e:W_e\text{ is infinite}\}$.
It is a standard fact that $\Tot\equiv_{\mathrm{tt}}\mathrm{INF}\equiv_{\mathrm{tt}}\emptyset''$. 

The following
fact explains the connection between the class $\BLR(A)$ and truth-table reducibility to the
jump and will be used throughout \Cref{sec:bounded-threshold}.

\begin{lemma}[Cole and Simpson~{\cite{ColeSimpson07}}]
\label{lem:tt-wtt-mindchange}
For sets $A$ and $B$, the following are equivalent.
\begin{enumerate}[label=\textnormal{(\roman*)}]
\item $A\leq_{\mathrm{tt}}B'$.
\item $A\leq_{\mathrm{wtt}}B'$.
\item The characteristic function of $A$ belongs to $\BLR(B)$.
\end{enumerate}
\end{lemma}

\begin{remark}
\label{rem:degree-invariance}
Cole and Simpson~\cite{ColeSimpson07} observed that
$A\leq_TB$ implies $\BLR(A)\subseteq\BLR(B)$. In particular,
$A\equiv_TB$ implies $\BLR(A)=\BLR(B)$. 
\end{remark}

\section{The robust properties}
\label{sec:density-degrees}

\begin{proposition}
\label{prop:coarse-hyperimmune}
\label{prop:positive-density-ioe}
For every Turing degree $\mathbf a$, the following are equivalent.
\begin{enumerate}[label=\textnormal{(\roman*)}]
\item $\mathbf a$ is hyperimmune.
\item $\mathbf a$ computes a
function that is not coarsely recursively bounded.
\item There is a function $h\leq_T\mathbf a$ such that
$\overline\rho(\{n:h(n)=g(n)\})=1$
for every recursive function $g$.
\end{enumerate}
\end{proposition}

\begin{proof}
\textnormal{(ii)$\Rightarrow$(i).}
Every eventual recursive bound is also a coarse recursive bound.

\textnormal{(i)$\Rightarrow$(iii).}
By the infinitely-often-equal characterization of hyperimmune degrees
(see Monin and Nies~\cite{MR3538473}), there is a function
$f\leq_T\mathbf a$ that agrees with every recursive function at infinitely
many arguments. Fix an effective coding of finite sequences with recursive
decoding, and use the blocks $\mathcal Q_k$ from
\eqref{eq:quadratic-blocks}. Set $h(0)=h(1)=0$.
If $f(k)$ codes a sequence of length $|\mathcal Q_k|$, use that sequence
as the values of $h$ on $\mathcal Q_k$, in increasing order of the inputs;
otherwise, set $h$ equal to zero throughout $\mathcal Q_k$.
This defines a total function $h\leq_T f$.

Fix a recursive function $g$. The function
$G(k)=\operatorname{code}\bigl(g(q_k),\ldots,g(q_{k+1}-1)\bigr)$
is recursive. By the choice of $f$, $f(k)=G(k)$ for infinitely many $k$.
For every such $k$, the functions $h$ and $g$ agree throughout
$\mathcal Q_k$. Hence \Cref{lem:long-block-density}\textnormal{(ii)} gives
$\overline\rho(\{n:h(n)=g(n)\})=1$.

\textnormal{(iii)$\Rightarrow$(ii)} is immediate.
\end{proof}

Reversing the direction of domination
gives an exact density analogue of Martin's characterization of the high
degrees.

\begin{proposition}
\label{thm:density-martin}
For every Turing degree $\mathbf a$, the following are equivalent.
\begin{enumerate}[label=\textnormal{(\roman*)}]
\item $\mathbf a$ is high.
\item There is a function $f\leq_T\mathbf a$ such that $g\ll_c f$ for every
recursive function $g$.
\item There is a function $f\leq_T\mathbf a$ such that
$\underline\rho(\{n:f(n)\geq g(n)\})>0$
for every recursive function $g$.
\end{enumerate}
\end{proposition}

\begin{proof}
Martin's domination theorem implies
\textnormal{(i)}$\Rightarrow$\textnormal{(ii)}. The implication
\textnormal{(ii)}$\Rightarrow$\textnormal{(iii)} is immediate. It remains to prove \textnormal{(iii)}$\Rightarrow$\textnormal{(i)}. Let
$f\leq_T\mathbf a$ witness \textnormal{(iii)}. Using the blocks
$\mathcal Q_n$ from \eqref{eq:quadratic-blocks}, define $F(n)=\max\{f(m):m\in\mathcal Q_n\}$.
Then $F\leq_T f\leq_T\mathbf a$. Fix a recursive function $g$.
\textnormal{(iii)} gives positive lower density to the set
$C=\{m:f(m)\geq g^{\mathrm{blk}}(m)\}$.
By \Cref{lem:long-block-density}\textnormal{(i)} with $\eta=0$, we have
$C\cap\mathcal Q_n\neq\varnothing$ for every sufficiently large~$n$.
For each such $n$, any $m\in C\cap\mathcal Q_n$ satisfies $F(n)\geq f(m)\geq g^{\mathrm{blk}}(m)=g(n)$.
Thus $F$ dominates every recursive function, and hence $\mathbf a$ is high.
\end{proof}

Following Downey, Jockusch and Stob \cite{DJS}, a set $A$ is
\emph{array recursive} if there is a
single function $g\leq_{\mathrm{wtt}}\emptyset'$ that dominates every
function $f\leq_TA$.  The same construction, with the density
step supplied by \Cref{lem:long-block-density}\textnormal{(i)}, shows that this
notion is unchanged by passing to asymptotic density.

\begin{proposition}
\label{thm:density-array-computable}
Let $A\subseteq\N$. The following
are equivalent.
\begin{enumerate}[label=\textnormal{(\roman*)}]
\item $A$ is array recursive.
\item There is a function $g\leq_{\mathrm{wtt}}\emptyset'$ such that
$f\ll_c g$ for every $f\leq_TA$.
\item There is a function $g\leq_{\mathrm{wtt}}\emptyset'$ such that
$\underline\rho(\{n:g(n)\geq f(n)\})>0$
for every $f\leq_TA$.
\end{enumerate}
\end{proposition}

Replacing lower by upper density in
\Cref{thm:density-array-computable} yields a much larger class.

\begin{proposition}
\label{prop:upper-density-array-incompleteness}
Let $A\subseteq\N$. The following are equivalent.
\begin{enumerate}[label=\textnormal{(\roman*)}]
\item $\emptyset'\nleq_TA$.
\item There is a function $g\leq_{\mathrm{wtt}}\emptyset'$ such that
$\overline\rho(\{n:g(n)\geq f(n)\})=1$
for every $f\leq_TA$.
\item There is a function $g\leq_{\mathrm{wtt}}\emptyset'$ such that
$ \overline\rho(\{n:g(n)\geq f(n)\})>0$
for every $f\leq_TA$.
\end{enumerate}
\end{proposition}

\begin{proof}
The implications
\textnormal{(ii)}$\Rightarrow$\textnormal{(iii)}$\Rightarrow$\textnormal{(i)} are immediate. It remains to prove
\textnormal{(i)}$\Rightarrow$\textnormal{(ii)}. Let $(K_s)_{s\in\N}$ be
the standard increasing recursive approximation to $K=\emptyset'$, and
define 
$m(n)=\min\{s:K_s\mathbin{\upharpoonright}n
                   =K\mathbin{\upharpoonright}n\}$.
Then $m\leq_{\mathrm{wtt}}K$ and no function $u\leq_TA$ dominates $m$ since $\emptyset'\nleq_TA$.
Put $g=m^{\mathrm{blk}}\leq_{\mathrm{wtt}}\emptyset'$. 

Fix $f\leq_TA$ and define
$ F(n)=\max\{f(k):k\in\mathcal Q_n\}$.
Since $F\leq_TA$, it does not dominate $m$. Hence there are infinitely
many $n$ such that $m(n)>F(n)$. For each such $n$,
$\mathcal Q_n\subseteq\{k:g(k)\geq f(k)\}$.
Consequently, \Cref{lem:long-block-density}\textnormal{(ii)} gives
$\overline\rho(\{k:g(k)\geq f(k)\})=1$, as required.
\end{proof}

The same lower-density mechanism also gives a partial-function form of
lowness.

\begin{proposition}
\label{prop:density-low}
Let $A\leq_T\emptyset'$. The following are equivalent.
\begin{enumerate}[label=\textnormal{(\roman*)}]
\item $A$ is low.
\item There is a function $f\leq_T\emptyset'$ such that, for every partial
$A$-recursive function $\psi$, $(\forall^\infty n)\,[\psi(n)\mathord\uparrow\ \text{or}\ \psi(n)<f(n)]$.
\item There is a function $f\leq_T\emptyset'$ such that, for every partial
$A$-recursive function $\psi$, $\rho\bigl(\{n:\psi(n)\mathord\uparrow\ \text{or}\ \psi(n)<f(n)\}\bigr)=1$.
\item There is a function $f\leq_T\emptyset'$ such that, for every partial
$A$-recursive function $\psi$,
$\underline\rho\bigl(
    \{n:\psi(n)\mathord\uparrow\ \text{or}\ \psi(n)<f(n)\}
  \bigr)>0$.
\end{enumerate}
\end{proposition}

\begin{proof}
It is well-known that \textnormal{(i)} and
\textnormal{(ii)} are equivalent, and
\textnormal{(ii)}$\Rightarrow$\textnormal{(iii)}$\Rightarrow$\textnormal{(iv)}
is immediate.
The implication \textnormal{(iv)}$\Rightarrow$\textnormal{(ii)} follows by
the same argument as
\textnormal{(iii)}$\Rightarrow$\textnormal{(i)} in
\Cref{thm:density-martin}, now applied to partial $A$-recursive functions.
\end{proof}

Replacing lower density by upper density in the preceding proposition
characterizes incompleteness.
For a binary tree $T\subseteq2^{<\N}$, closed under taking initial
segments, its \emph{width} is $\sup_n|T\cap2^n|$; thus $T$ has
\emph{bounded width} if a single finite constant bounds the number of
nodes on every level. We write $[T]$ for the set of its infinite paths.

\begin{lemma}[Kaufmann and
Kummer~\cite{KaufmannKummer96}]
\label{lem:bounded-width-paths}
Let $A\subseteq\N$. If $T\subseteq2^{<\N}$ is an $A$-r.e. binary tree
of bounded width, then every $X\in[T]$ is $A$-recursive.
\end{lemma}

\begin{theorem}
\label{prop:upper-density-partial-incompleteness}
Let $A\leq_T\emptyset'$. The following are equivalent.
\begin{enumerate}[label=\textnormal{(\roman*)}]
\item $A<_T\emptyset'$.
\item There is a function $f\leq_T\emptyset'$ such that, for every partial
function $\psi$ recursive in $A$,
$\overline\rho\bigl(
    \{n:\psi(n)\mathord\uparrow\ \text{or}\ \psi(n)<f(n)\}
  \bigr)=1$.
\item There is a function $f\leq_T\emptyset'$ such that, for every partial
function $\psi$ recursive in $A$,
$\overline\rho\bigl(
    \{n:\psi(n)\mathord\uparrow\ \text{or}\ \psi(n)<f(n)\}
  \bigr)>0$.
\end{enumerate}
\end{theorem}

\begin{proof}
The implications
\textnormal{(ii)}$\Rightarrow$\textnormal{(iii)}$\Rightarrow$\textnormal{(i)} are immediate. 
Conversely, suppose that $A<_T\emptyset'$.  Let $(K_s)_{s\in\N}$ be the
standard increasing recursive approximation to $K=\emptyset'$, and let
$ m(n)=\min\{s:K_s\mathbin{\upharpoonright}n
                   =K\mathbin{\upharpoonright}n\}$.  Put $f=m^{\mathrm{blk}}\leq_{\mathrm{wtt}}\emptyset'$.
Suppose, toward a contradiction, that some partial $A$-recursive function
$\psi$ has a failure set
$E=\{k:\psi(k)\mathord\downarrow\ \text{and}\ \psi(k)\geq f(k)\}$
of positive lower density $\delta=\underline\rho(E)>0$.
By \Cref{lem:long-block-density}\textnormal{(i)}, for all sufficiently large
$n$,
$|E\cap\mathcal Q_n|
  >\frac{\delta}{2}|\mathcal Q_n|$.
Choose an integer $r$ so large that
$1/r<\delta/4$, and put
$t_n=\left\lceil\frac{|\mathcal Q_n|}{r}\right\rceil$.
Since $|\mathcal Q_n|\to\infty$,
$t_n <\frac{\delta}{2}|\mathcal Q_n|$ for all sufficiently large~$n$.

For each $n$, enumerate into $V_n\subseteq 2^n$ a string $\sigma$ as soon as
there are at least $t_n$ distinct $k\in\mathcal Q_n$ such that $\psi(k)\mathord\downarrow$ and $K_{\psi(k)}\mathbin{\upharpoonright}n=\sigma$.
The sequence $(V_n)_{n\in\N}$ is uniformly $A$-r.e., and
$|V_n|\leq r$ for every $n$, since the sets of inputs counted for distinct
strings are disjoint.

If $k\in E\cap\mathcal Q_n$, then
$\psi(k)\geq f(k)=m(n)$,
and therefore
$K_{\psi(k)}\mathbin{\upharpoonright}n
  =K\mathbin{\upharpoonright}n$.
The preceding lower bound on $|E\cap\mathcal Q_n|$ consequently implies
$K\mathbin{\upharpoonright}n\in V_n$
for every sufficiently large $n$.  Replacing finitely many initial
components by the corresponding singletons, we may assume that this holds
for every $n$.

Now form the $A$-r.e. tree
$T=\{\sigma\in 2^{<\N}:
       (\forall j\leq|\sigma|)\,
       \sigma\mathbin{\upharpoonright}j\in V_j\}$.
The tree $T$ has width at most $r$, and $K\in[T]$. By
\Cref{lem:bounded-width-paths}, $K\leq_TA$, a contradiction. Hence \textnormal{(i)}$\Rightarrow$\textnormal{(ii)}.
\end{proof}

We next consider agreement with recursive functions. A set $A$ is
\emph{recursively semi-traceable} if, for every $f\leq_TA$, there is a
recursive function $g$ such that $f(n)=g(n)$ for infinitely many~$n$; see
\cite{KHNSY10}.

\begin{proposition}
\label{thm:density-recursive-semitraceability}
Let $A\subseteq\N$. The following
are equivalent.
\begin{enumerate}[label=\textnormal{(\roman*)}]
\item $A$ is recursively semi-traceable.
\item For every $f\leq_TA$, there is a recursive function $g$ such that
$\overline\rho(\{n:f(n)=g(n)\})=1$.
\item For every $f\leq_TA$, there is a recursive function $g$ such that
$\overline\rho(\{n:f(n)=g(n)\})>0$.
\end{enumerate}
\end{proposition}

\begin{proof}
It suffices to prove \textnormal{(i)}$\Rightarrow$\textnormal{(ii)}. Given
$f\leq_TA$, let $F(n)$ code the finite sequence
$f\upharpoonright q_n$, where $(q_n)$ is fixed in
\eqref{eq:quadratic-blocks}. By recursive semi-traceability,
some recursive function $\widehat F$ agrees with $F$ at infinitely many
arguments. Define a recursive function $g$ as follows. If
$n\geq1$ and $m\in\mathcal Q_{n-1}$, decode $\widehat F(n)$; if it codes a
sequence $\sigma$ of length $q_n$, let $g(m)=\sigma(m)$, and otherwise let
$g(m)=0$. Give $g$ arbitrary values below~$q_0$.
Whenever $n\geq1$ and $\widehat F(n)=F(n)$, the functions $f$ and $g$ agree throughout
$\mathcal Q_{n-1}$.  Since this happens for infinitely many~$n$,
\Cref{lem:long-block-density}\textnormal{(ii)} shows that the agreement set
has upper density one.
\end{proof}

Replacing upper density by lower density forces recursiveness.

\begin{proposition}
\label{prop:ld-recursive}
A set $A\subseteq\N$ is recursive if and only if, for every $f\leq_TA$,
there is a recursive function $g$ such that
$\underline\rho(\{n:f(n)=g(n)\})>0$.
\end{proposition}

\begin{proof}
The forward implication is immediate. For the converse, suppose that $A$
satisfies the stated agreement condition. Fix a recursive coding
$\sigma\mapsto\ulcorner\sigma\urcorner$ of the finite binary strings by
natural numbers. Let $p(n)=\ulcorner A\upharpoonright n\urcorner$ and put
$F=p^{\mathrm{blk}}$. By assumption,
there is a recursive
function $g$ such that, for
$C=\{m:F(m)=g(m)\}$,
we have $\delta:=\underline\rho(C)>0$. Choose an integer $K\geq1$ such that
$1/K<\delta$.  By
\Cref{lem:long-block-density}\textnormal{(i)}, for all sufficiently large
$n$, $|C\cap\mathcal Q_n|
  >\frac{|\mathcal Q_n|}{K}$.
For each $n$, let
\[
  S_n=\left\{\sigma\in2^n:
  \bigl|\{m\in\mathcal Q_n:g(m)=\ulcorner\sigma\urcorner\}\bigr|
  >\frac{|\mathcal Q_n|}{K}\right\}.
\]
The sequence $(S_n)_{n\in\N}$ is uniformly recursive, and $|S_n|<K$ for
every $n$, since the sets of positions counted for distinct strings are
disjoint. Moreover, $A\upharpoonright n\in S_n$
for all sufficiently large $n$.
Choose $n_0$ such that $A\upharpoonright n\in S_n$ holds for every
$n\geq n_0$, and define a recursive binary tree
\[
  \mathcal P=\left\{\tau\in2^{<\N}:
  (\forall k)\bigl(n_0\leq k\leq|\tau|
  \Longrightarrow \tau\upharpoonright k\in S_k\bigr)\right\}.
\]
Then $A$ is an infinite path through $\mathcal P$, while
$|\mathcal P\cap2^n|<K$
for every $n\geq n_0$. Thus $\mathcal P$ has bounded width, so
$A$ is recursive by \Cref{lem:bounded-width-paths}.
\end{proof}

We now turn to strong jump traceability. An \emph{r.e. trace} is a uniformly
r.e. sequence $\mathcal T=(T_n)_{n\in\N}$ of finite subsets of $\N$.
An \emph{order function} is a positive, nondecreasing, unbounded recursive
function. A trace is \emph{$h$-bounded} if $|T_n|\leq h(n)$ for every~$n$,
and it traces a partial function $\psi$ if $\psi(n)\in T_n$ whenever
$\psi(n)\mathord\downarrow$. For such a function and trace, put
\begin{align*}
  \operatorname{Hit}(\psi,\mathcal T)
  &=\{n\in\operatorname{dom}(\psi):\psi(n)\in T_n\},\\
  \operatorname{Corr}(\psi,\mathcal T)
  &=\{n:\psi(n)\mathord\uparrow
       \text{ or }\psi(n)\in T_n\}.
\end{align*}
Thus $\mathcal T$ traces $\psi$ precisely when
$\operatorname{Corr}(\psi,\mathcal T)=\N$.  Recall that $A$ is
\emph{strongly jump traceable} if for every order function
$h$, every partial $A$-recursive function has an $h$-bounded r.e. trace; see
\cite{FigueiraNiesStephan08}.

\begin{theorem}
\label{thm:density-strong-jump-traceability}
A set $A$ is strongly jump traceable if and only if for every order function $h$ and every partial $A$-recursive function
$\psi$, there is an $h$-bounded r.e. trace $\mathcal T$ such that $\underline\rho(\operatorname{Corr}(\psi,\mathcal T))>0$.
\end{theorem}

\begin{proof}
The forward implication is immediate. Conversely, fix an order function $g$
and a partial $A$-recursive function $\psi$. To respect the prescribed bound
$g$, we use a slower auxiliary order and retain only values that occur
frequently within each block. Let
$b(n)=\max\{1,\lfloor\sqrt{g(n)}\rfloor\}$,
and put $a=b^{\mathrm{blk}}$.
Applying the assumed density condition to $a$ and $\psi^{\mathrm{blk}}$, there
is an $a$-bounded r.e. trace
$\mathcal V=(V_m)_{m\in\N}$ such that
$C=\operatorname{Corr}(\psi^{\mathrm{blk}},\mathcal V)$
has positive lower density.  Put $\delta=\underline\rho(C)>0$.  By
\Cref{lem:long-block-density}\textnormal{(i)}, for all sufficiently large
$n$,
$|C\cap\mathcal Q_n|
  >\frac{\delta}{2}|\mathcal Q_n|$.
If $n\in\operatorname{dom}(\psi)$, the members of
$C\cap\mathcal Q_n$ are precisely the $m\in\mathcal Q_n$ for which
$\psi(n)\in V_m$.
Let
$k_n=\left\lceil
       \frac{|\mathcal Q_n|b(n)}{g(n)}
      \right\rceil$
and enumerate into $S_n$ each value $y$ that appears in at least $k_n$
distinct components $V_m$ with $m\in\mathcal Q_n$.  The sequence
$(S_n)_{n\in\N}$ is uniformly r.e., and
\[
  |S_n|k_n
  \leq\sum_{m\in\mathcal Q_n}|V_m|
  \leq |\mathcal Q_n|b(n)
  \leq g(n)k_n.
\]
Consequently, $|S_n|\leq g(n)$.  Moreover,
$\frac{k_n}{|\mathcal Q_n|}
  \leq\frac{b(n)}{g(n)}+\frac1{|\mathcal Q_n|}
  \longrightarrow0$.
It follows from the preceding lower bound on
$|C\cap\mathcal Q_n|$ that $\psi(n)\in S_n$
for almost every $n\in\operatorname{dom}(\psi)$.  After replacing finitely
many components by singletons containing the corresponding convergent
values, $(S_n)_{n\in\N}$ becomes a $g$-bounded r.e. trace for $\psi$.
Since $g$ and $\psi$ were arbitrary, $A$ is strongly jump traceable.
\end{proof}

To formulate an upper-density analogue, we first need to account for
the possible sparsity of the domain. Since
$\operatorname{Hit}(\psi,\mathcal T)\subseteq\operatorname{dom}(\psi)$,
an infinite domain of density zero prevents even a trace that captures
every defined value from having positive upper density of hits.
We therefore pass to $\psi^{\mathrm{blk}}$: its domain contains
$\mathcal Q_n$ for every $n\in\operatorname{dom}(\psi)$, and hence has
upper density one whenever $\operatorname{dom}(\psi)$ is infinite, by
\Cref{lem:long-block-density}\textnormal{(ii)}.
We must also use $\operatorname{Hit}$ rather than
$\operatorname{Corr}$, so that undefined values do not count as successes.
Recall that a total function $d$ is \emph{diagonally non-recursive}, or
\emph{DNR}, if $d(e)\neq\Phi_e(e)$ whenever
$\Phi_e(e)\mathord\downarrow$. A degree is \emph{non-DNR} if it computes
no DNR function.
We have the following characterization due to Kjos-Hanssen, Merkle, and
Stephan~\cite{KHMS11} (see H\"olzl and
Merkle~\cite[Theorems~10 and~20]{HoelzlMerkle10}). A set $A$ computes
no DNR function if and only if there is an order function $b$ such that
every total function $f\leq_TA$ has a $b$-bounded r.e. trace
$(V_n)_{n\in\N}$ with $f(n)\in V_n$ for infinitely many~$n$.
Moreover, if $A$ computes no DNR function, then for every order function
$b$ and every partial $A$-recursive function $\psi$ with infinite domain,
there is a $b$-bounded r.e. trace $(V_n)_{n\in\N}$ such that
$\psi(n)\in V_n$ for infinitely many $n\in\operatorname{dom}(\psi)$.

The following proposition follows readily from these results by a block
argument.

\begin{proposition}
\label{prop:upper-density-strong-jump-traceability}
Let $A\subseteq\N$. The following are equivalent.
\begin{enumerate}[label=\textnormal{(\roman*)}]
\item $A$ computes no DNR function.
\item For every order function $h$ and every partial $A$-recursive function
$\psi$ with infinite domain, there is an $h$-bounded r.e. trace
$\mathcal T$ such that
$\overline\rho(\operatorname{Hit}(\psi^{\mathrm{blk}},\mathcal T))=1$.
\item For every order function $h$ and every partial $A$-recursive function
$\psi$ with infinite domain, there is an $h$-bounded r.e. trace
$\mathcal T$ such that
$\overline\rho(\operatorname{Hit}(\psi^{\mathrm{blk}},\mathcal T))>0$.
\end{enumerate}
\end{proposition}

Let $\mu$ denote the usual product measure on $2^\N$. Recall that a set
$A$ is \emph{almost-everywhere dominating} (\emph{a.e.d.}) if, for $\mu$-almost every
$X\in2^\N$ and every total function $g\leq_T X$, there is a function
$f\leq_TA$ which dominates $g$. It is \emph{uniformly almost-everywhere
dominating} (\emph{u.a.e.d.}) if one function $f\leq_TA$ works for almost every $X$ and all
total $g\leq_T X$; see Dobrinen and
Simpson~\cite{DobrinenSimpson04}.
For a Turing functional $\Phi$, let
$\Tot(\Phi)=\{X\in2^\N:\Phi^X\text{ is total}\}$.
If $f$ is total, define the domination set
\begin{equation*}
  D_{\Phi,f}(X)
  =\{n:\Phi^X(n)\mathord\downarrow\text{ and }\Phi^X(n)\leq f(n)\}.
\end{equation*}

\begin{definition}[Kjos-Hanssen~\cite{KjosHanssen07}]
\label{def:pmd}
A set $A$ is \emph{positive-measure dominating} (\emph{p.m.d.}) if, for every $\Phi$
with $\mu(\Tot(\Phi))>0$, there is an $f\leq_TA$ such that
$\mu\bigl(\{X\in\Tot(\Phi):f\text{ dominates }\Phi^X\}\bigr)>0$.
\end{definition}

\begin{theorem}[Kjos-Hanssen, Miller, and Solomon
{\cite{KMS12}}]
\label{thm:kms-domination}
For every set $A$, the following are equivalent.
\begin{enumerate}[label=\textnormal{(\roman*)}]
\item $A$ is u.a.e.d.
\item $A$ is a.e.d.
\item $A$ is p.m.d.
\end{enumerate}
\end{theorem}

We now replace eventual domination in the positive-measure definition by a
density condition.

\begin{proposition}
\label{thm:aed-density}
For every set $A$, the following are equivalent.
\begin{enumerate}[label=\textnormal{(\roman*)}]
\item $A$ is a.e.d.
\item For every Turing functional $\Phi$ with
$\mu(\Tot(\Phi))>0$, there is an $f\leq_TA$ such that
$\mu\bigl(\{X\in\Tot(\Phi):
       \underline\rho(D_{\Phi,f}(X))>0\}\bigr)>0$.
\item There is an $f\leq_TA$ such that
\[
 \mu\bigl(\{X:(\forall g\leq_T X)\;
   \underline\rho(\{n:g(n)\leq f(n)\})>0\}\bigr)=1.
\]
\end{enumerate}
\end{proposition}

\begin{proof}
\textnormal{(i)$\Rightarrow$(iii).}
By \Cref{thm:kms-domination}, $A$ is u.a.e.d. and hence (iii) holds.

\textnormal{(iii)$\Rightarrow$(ii)} is immediate.

\textnormal{(ii)$\Rightarrow$(i).}
Assume \textnormal{(ii)} and fix $\Phi$ with
$\mu(\Tot(\Phi))>0$. Let
$\widehat\Phi^X(m)=\max_{i\leq m}\Phi^X(i)$,
with the value defined only when all computations on the right converge.
Then $\Tot(\widehat\Phi)=\Tot(\Phi)$. Applying the density condition to
$\widehat\Phi$, choose $f\leq_TA$ and a positive-measure set
$P\subseteq\Tot(\Phi)$ such that
$C_X=\{m:f(m)\geq\widehat\Phi^X(m)\}$
has positive lower density for every $X\in P$.
Define $F\leq_TA$ by
$F(n)=\max\{f(m):m\in\mathcal Q_n\}$.
By \Cref{lem:long-block-density}\textnormal{(i)}, for every $X\in P$, the set
$C_X$ meets $\mathcal Q_n$ for all
sufficiently large $n$.
Choosing $m\in C_X\cap\mathcal Q_n$ and using $m\geq q_n\geq n$ gives
$F(n)\geq f(m)\geq\widehat\Phi^X(m)\geq\Phi^X(n)$.
Hence $F$ witnesses p.m.d. for $\Phi$. Since $\Phi$ was arbitrary,
$A$ is p.m.d., therefore a.e.d. by \Cref{thm:kms-domination}.
\end{proof}

Replacing lower density by upper density in
\Cref{thm:aed-density}\textnormal{(iii)} gives
hyperimmunity.

\begin{proposition}
\label{prop:upper-density-aed-hyperimmune}
For every set $A$, the following are equivalent.
\begin{enumerate}[label=\textnormal{(\roman*)}]
\item $A$ has hyperimmune degree.
\item There is an $f\leq_TA$ such that
\[
 \mu\bigl(\{X:(\forall g\leq_T X)\;
   \overline\rho(\{n:g(n)\leq f(n)\})=1\}\bigr)=1.
\]
\item There is an $f\leq_TA$ such that
\[
 \mu\bigl(\{X:(\forall g\leq_T X)\;
   \overline\rho(\{n:g(n)\leq f(n)\})>0\}\bigr)=1.
\]
\end{enumerate}
\end{proposition}

\begin{proof}
Suppose first that $A$ has hyperimmune degree, and choose a hyperimmune
function $h\leq_TA$.  For almost every $X$, the function $h$ remains
$X$-hyperimmune.  Indeed, if for some functional $\Phi$ and some $k$ the
class of $X$ satisfying $\Phi^X(n)\mathord\downarrow\geq h(n)$ for every
$n\geq k$ had positive measure, then, inside a cylinder on which this class
has relative measure greater than $3/4$, one could compute an eventual
bound for $h$, a
contradiction. Taking the countable union over all functionals $\Phi$ and
all $k$ proves the assertion.
Put $f=h^{\mathrm{blk}}$.  For every such
$X$ and every $g\leq_T X$, $h(n)>\max\{g(m):m\in\mathcal Q_n\}$ for infinitely many $n$, and
for each such $n$ the whole block $\mathcal Q_n$ is contained in
$\{m:g(m)\leq f(m)\}$.  By
\Cref{lem:long-block-density}\textnormal{(ii)}, this set has upper density
one.  Thus \textnormal{(i)} implies \textnormal{(ii)}, while
\textnormal{(ii)} implies \textnormal{(iii)} and \textnormal{(iii)} implies \textnormal{(i)} trivially.
\end{proof}

The methods used in this section apply more generally. Similar arguments
yield density versions of other notions, including recursive traceability,
r.e. traceability, and jump traceability. We omit the details.

\section{Bounded-threshold domination and superhighness}
\label{sec:bounded-threshold}

We now answer the question of Kjos-Hanssen and Nies~\cite{KHN09}.

\begin{theorem}[with Liang Yu]
\label{thm:superhigh-threshold}
\label{cor:clean-threshold}
A set $A$ is superhigh if and only if there are functions
$h\in\BLR(A)$ and $f\leq_TA$ such that
\begin{equation}\label{eq:blr-threshold-domination}
  (\forall e\in\Tot)(\forall n\geq h(e))\;
  \Phi_e(n)\leq f(n).
\end{equation}
\end{theorem}

\begin{proof}
Suppose first that $h\in\BLR(A)$ and $f\leq_TA$ satisfy
\eqref{eq:blr-threshold-domination}. Choose an $A$-recursive approximation
$H(e,s)$ to $h(e)$ and a recursive function $g$ such that
$|\{s:H(e,s)\neq H(e,s+1)\}|<g(e)$
for every $e$. Replacing $H$ by
$\widehat H(e,s)=\max_{t\leq s}H(e,t)$,
and $h$ by $\widehat h(e)=\lim_s\widehat H(e,s)$, we may assume without
loss of generality that $H(e,s)\leq H(e,s+1)$ for all $e,s$. Indeed,
$\widehat h\in\BLR(A)$ and $\widehat h(e)\geq h(e)$, so
\eqref{eq:blr-threshold-domination} remains true.

Let $\theta$ be a recursive function such that $\Phi_{\theta(e)}(n)$ is the
least stage $t$, if it exists, for which $\Phi_e(m)[t]\mathord\downarrow$ for all $m\leq n$.
Thus, if $\Phi_e$ is total, then $\Phi_{\theta(e)}$ is total, and $\Phi_e(m)[\Phi_{\theta(e)}(n)]\mathord\downarrow$ for all $m \le n$.

We define an $A$-recursive binary approximation $T(n,s)$ to $\Tot$. Put
$T(n,0)=0$. For $s>0$, let
$k_s=H(\theta(n),s)$.
If either $k_s\geq s$ or $k_s>H(\theta(n),s-1)$, put $T(n,s)=0$.
Otherwise, put $T(n,s)=1$ exactly when
\[
  (\forall m\leq k_s)\;\Phi_n(m)[f(k_s)]\mathord\downarrow
\]
and
\[
  (\forall m)\bigl(k_s<m<s\Longrightarrow
  \Phi_n(m)[f(m)]\mathord\downarrow\bigr);
\]
if either condition fails, put $T(n,s)=0$.

We first check that $\lim_sT(n,s)=\Tot(n)$. If $\Phi_n$ is total, then
$\Phi_{\theta(n)}$ is total. Let $k=h(\theta(n))$. By
\eqref{eq:blr-threshold-domination} applied to $\theta(n)$, for every
$m\geq k$ we have $\Phi_{\theta(n)}(m)\leq f(m)$.
Hence, for all sufficiently large $s$, after $H(\theta(n),\cdot)$ has
stabilized at $k$, both displayed conditions hold, and therefore
$T(n,s)=1$. If $\Phi_n$ is not total, then it is clear that for all sufficiently large stages $s$, one of the two displayed conditions fails. And therefore $T(n,s)=0$.

Finally, $T$ has a recursive mind-change bound. If
$H(\theta(n),\cdot)$ changes $c$ times, then $c<g(\theta(n))$.
On each interval on which $H(\theta(n),s)$ is constant, say with value $k$,
the sequence $T(n,s)$ has the form $0^*1^*0^*$: the reset rule and the
condition $k\geq s$ initially force the value $0$, and thereafter the tests
only acquire additional conjuncts as $s$ increases. Thus there are at most
two changes within each of the $c+1$ constant intervals, and at most one
additional change at each of their $c$ boundaries. Consequently,
\[
  |\{s:T(n,s)\neq T(n,s+1)\}|
  \leq 2(c+1)+c=3c+2<3g(\theta(n))+3.
\]
By \Cref{lem:tt-wtt-mindchange}, $\Tot\leq_{\mathrm{tt}}A'$ and $A$ is superhigh.

Conversely, suppose that $A$ is superhigh. Fix a
functional $\Gamma$ and a recursive function $\lambda$ such
that, for every $e$,
$\Tot(e)=\Gamma^{A'\upharpoonright\lambda(e)}(e)$
and $\Gamma^\sigma(e)\in \{0,1\}$ for every finite string
$\sigma$ of length $\lambda(e)$.

Let $A'_s$ be the standard $A$-recursive approximation to $A'$. We define
$H\leq_TA$ and $f\leq_TA$ in stages. At stage $0$, put $H(i,0)=i$ for
every $i$.

At stage $s>0$, first update $H$. If there is no $e<s$ such that $A'_s\upharpoonright\lambda(e)
  \neq A'_{s-1}\upharpoonright\lambda(e)$,
put $H(i,s)=H(i,s-1)$ for all $i$. Otherwise, let $e_s$ be the least such
$e$, and define
\begin{equation}\label{eq:threshold-stage-update}
  H(i,s)=
  \begin{cases}
    H(i,s-1),&\text{if }i<e_s,\\
    \max\{H(i,s-1),s+i\},&\text{if }i\geq e_s.
  \end{cases}
\end{equation}
Thus $H(i,\cdot)$ is nondecreasing for every $i$.

We next define $f$. At stage $s$, let $n$ be the least number such that
$f(n)$ has not yet been defined. Let
\begin{equation*}
  E_{n,s}=\{e<s:H(e,s)\leq n
  \text{ and }\Gamma^{A'_s\upharpoonright\lambda(e)}(e)=1\}.
\end{equation*}
If $\Phi_e(n)[s]\mathord\downarrow$ for every $e\in E_{n,s}$, define
\begin{equation*}
  f(n)=\max\bigl(\{0\}\cup
  \{\Phi_e(n)[s]:e\in E_{n,s}\}\bigr).
\end{equation*}
Otherwise define no new value of $f$ at stage $s$.

We verify that the construction works. First, $h(i)=\lim_sH(i,s)$ exists for
every $i$. Indeed, $H(i,\cdot)$ can change only at stages at which, for
some $e\leq i$, the finite string
$A'_s\upharpoonright\lambda(e)$ changes. For each such $e$, this finite
string changes at most $\lambda(e)$ times. Thus
$|\{s:H(i,s)\neq H(i,s+1)\}|
  \leq\sum_{e\leq i}\lambda(e)$. Since the
construction of $H$ uses only the oracle $A$, it follows that
$h\in\BLR(A)$.

Next we show that $f$ is total. We argue by induction on $n$. Assume that
$f\upharpoonright n$ has already been defined. Since $H(i,0)=i$ and
$H(i,\cdot)$ is nondecreasing, if $H(i,s)\leq n$ for some $s$, then
$i\leq n$. 
Choose a stage $s_0>n$ after which all strings
$A'_s\upharpoonright\lambda(i)$ with $i\leq n$ and all values $H(i,s)$
with $i\leq n$ have stabilized. Then, for all $s\geq s_0$, the set
$E_{n,s}$ is a fixed finite set consisting only of indices $e$ for which
$\Tot(e)=1$. Thus every $\Phi_e$ with $e\in E_{n,s}$ is total, and there
is a later stage $s_1$ such that
$\Phi_e(n)[s_1]\mathord\downarrow$ for every $e$ in this finite set. If
$f(n)$ has not already been defined, the construction defines it at stage
$s_1$. Therefore $f$ is total and
$f\leq_TA$.

It remains to verify \eqref{eq:blr-threshold-domination}. Fix
$e\in\Tot$ and $n\geq h(e)$. Let $t$ be the stage at which $f(n)$ is
defined. By construction, $t>n$. Since $H(e,\cdot)$ is
nondecreasing and has limit $h(e)$,
$H(e,t)\leq h(e)\leq n$.
Also $H(e,0)=e$, so $e\leq h(e)\leq n<t$.

We claim that
$A'_t\upharpoonright\lambda(e)=A'\upharpoonright\lambda(e)$. Otherwise,
this string would change at some later stage $u>t$. At that stage the least
changed index $e_u$ would be at most $e$, and the update rule
\eqref{eq:threshold-stage-update} would give
$H(e,u)\geq u+e>n$,
contradicting the fact that the nondecreasing sequence $H(e,\cdot)$ has
final value at most $n$.

Since $\Phi_e$ is total, $e\in E_{n,t}$. When $f(n)$ was
defined at stage $t$, the construction required
$\Phi_e(n)[t]\mathord\downarrow$ for every $e\in E_{n,t}$ and 
$\Phi_e(n)\leq f(n)$. Since
$n\geq h(e)$ was arbitrary, \eqref{eq:blr-threshold-domination} holds.
\end{proof}

The bounded-threshold characterization also gives an effective
strengthening of maximality. Martin~\cite{MR224469} characterized the
high r.e. degrees as exactly those containing maximal sets.
An r.e. coinfinite set $A$ is \emph{maximal} if its complement is
cohesive: for every $e$, either $\overline A\cap W_e$ or
$\overline A\setminus W_e$ is finite. We obtain a stronger condition by
requiring a bounded-limit threshold for this eventual homogeneity. Call
an r.e. coinfinite set $A$ \emph{supermaximal} if there is a function
$f\in\BLR(A)$ such that
\begin{equation}
\label{eq:supermaximal-threshold}
  (\forall e)(\forall n>m>f(e))\;
  \bigl[n,m\notin A\ \Longrightarrow
  (n\in W_e\Longleftrightarrow m\in W_e)\bigr].
\end{equation}
Thus every supermaximal set is maximal.

\begin{theorem}
\label{thm:superhigh-supermaximal}
Every supermaximal set is superhigh, and every superhigh r.e. degree
contains a supermaximal set.
\end{theorem}

\begin{proof}
First suppose that $A$ is supermaximal, witnessed by $f\in\BLR(A)$.
Choose a total $A$-recursive approximation $F(e,s)\to f(e)$ whose
number of changes is bounded by a recursive function $r(e)$.
Fix a recursive function $\theta$ such that
\[
  W_{\theta(e)}=
  \begin{cases}
    [0,|W_e|),&\text{if $W_e$ is finite},\\
    \N,&\text{if $W_e$ is infinite}.
  \end{cases}
\]
Using $A$, define $x(e,s)=\min\{x>F(\theta(e),s):x\notin A\}$. Also define $T(e,s)=\mathbf1_{W_{\theta(e),s}}(x(e,s))$.
These functions are total because $A$ is coinfinite. On each interval
on which $F(\theta(e),s)$ is constant, $x(e,s)$ is constant and
$T(e,s)$ changes at most once. Counting changes between such
intervals, $T(e,s)$ changes at most $2r(\theta(e))+1$ times.
Define $x_e=\min\{x>f(\theta(e)):x\notin A\}$.
If $W_e$ is infinite, then $T(e,s)\to1$. If $W_e$ is finite, then
$x_e\notin W_{\theta(e)}$: otherwise
\eqref{eq:supermaximal-threshold} would put every element of
$\overline A$ greater than $x_e$ into the finite set $W_{\theta(e)}$.
Thus $T(e,s)\to0$ in this case. Consequently,
$\mathbf1_{\mathrm{INF}}\in\BLR(A)$. By
\Cref{lem:tt-wtt-mindchange},
$\mathrm{INF}
  \leq_{\mathrm{tt}}A'$,
so $A$ is superhigh.

For the converse, let $D$ be an r.e. set of the prescribed superhigh
degree. By \Cref{thm:superhigh-threshold}, choose $d\leq_TD$ and
$h\in\BLR(D)$ such that
\begin{equation}
\label{eq:supermaximal-domination}
  (\forall j\in\Tot)(\forall x\geq h(j))\;
  \Phi_j(x)\leq d(x).
\end{equation}
We also use
$\BLR(D)=\BLR(\emptyset')$:
since $D\leq_T\emptyset'$ is superhigh, Cole and
Simpson~\cite[Theorem~6.12]{ColeSimpson07}, applied with
$X=\emptyset'$ and $Y=D$, give
$\BLR(\emptyset')\subseteq\BLR(D)$. The reverse inclusion follows
from $D\leq_T\emptyset'$.
Fix a recursive approximation $d_s(x)\to d(x)$ with a
$D$-recursive convergence modulus, see \cite{Soare87}.

We construct an r.e. set $A$ by defining an increasing sequence of
finite sets $(A_s)_{s\in\N}$, with $A=\bigcup_s A_s$. Our aim is to
make all sufficiently large elements of $\overline A$ agree on
membership in each $W_e$, with a threshold in $\BLR(A)$. To organize
these requirements, define the \emph{$i$-state} of $x$ at stage $s$ by
$\nu_i(x,s)=\sum_{j\leq i}2^{i-j}\mathbf1_{W_{j,s}}(x)$.
This is the binary code of the membership vector
$(\mathbf1_{W_{0,s}}(x),\ldots,\mathbf1_{W_{i,s}}(x))$.
Numerical comparison of states therefore orders these vectors
lexicographically, giving priority to the smaller indices. For each
fixed $i$ and $x$, the state is nondecreasing in $s$ and takes values
in $\{0,\ldots,2^{i+1}-1\}$. Moreover, for $e\leq i$,
$\nu_e(x,s)=\left\lfloor\frac{\nu_i(x,s)}{2^{i-e}}\right\rfloor$.
Thus replacing a marker by one with at least as large an $i$-state
does not decrease any of its $e$-states for $e\leq i$. We use this
ordering to select the markers that remain outside $A$.

At any point in stage $s$, the \emph{current markers} are the numbers
in $[0,s]$ that have not yet been enumerated into $A$. Write their
increasing list as $M=\{m_0<\cdots<m_{k-1}\}$;
the \emph{rank} of $m_i$ is its current index $i$ in this list.
A \emph{tail action at rank $i<k$} is defined as follows. Let
\[
  \ell=\min\bigl\{j:i\leq j<k\text{ and }
  \nu_i(m_j,s)=\max_{i\leq t<k}\nu_i(m_t,s)\bigr\}.
\]
Enumerate every element of $\{m_i,\ldots,m_{k-1}\}\setminus\{m_\ell\}$
into $A$.
After each action, ranks are reassigned according to the new
increasing list; all subsequent actions use these updated ranks.

Put $A_{-1}=D_{-1}=\emptyset$ and $d_{-1}(x)=0$ for every $x$.
At stage $s$, start with the marker list
$M=[0,s]\setminus A_{s-1}$ and perform the following steps.
\begin{enumerate}[label=\textnormal{(\arabic*)}]
\item \emph{Coding actions.}
Process the elements $n\in D_s\setminus D_{s-1}$ in increasing order.
For each such $n$, perform a tail action at rank $n$ if that rank
is currently present; otherwise do nothing.
\item \emph{Ordinary action.}
Using the current list $m_0<\cdots<m_{k-1}$ after all coding
actions, declare rank $i<k$ \emph{eligible} if there is a $j$ with
$i<j<k$ such that $\nu_i(m_j,s)>\nu_i(m_i,s)$, and $s\leq d_s(m_i) \lor d_s(m_i)\neq d_{s-1}(m_i) $.
If any rank is eligible, perform one tail action at the least
eligible rank. Otherwise perform no ordinary action.
\end{enumerate}
This completes stage $s$.

\begin{claim}
There are only finitely many actions at each rank, and $A$ is
coinfinite. The rank of any marker that survives an action cannot
increase.
\end{claim}

\begin{proof}[Proof of the claim]
We argue inductively on $i$. After the last action at a smaller rank
and the last coding action at rank $i$, every action at rank
$i$ strictly increases the $i$-state of its marker. Between actions,
that state cannot decrease. Since it takes at most $2^{i+1}$
possible values, there are only finitely many such actions. Once
the smaller ranks have stabilized, new numbers ensure that rank $i$
is eventually present and stabilizes as well. Hence $A$ is coinfinite.
Finally, an action only removes markers, so the rank of a surviving
marker cannot increase.
\end{proof}

Write
$\overline A=\{c_0<c_1<c_2<\cdots\}$.

\begin{claim}
\label{clm:supermaximal-degree}
$A\equiv_TD$.
\end{claim}

\begin{proof}[Proof of the claim]
For $A\leq_TD$, fix $z$. An ordinary action
that enumerates $z$ is triggered by a marker $x\leq z$. Using $D$,
find a stage beyond the convergence stages and final values of all
$d_s(x)$ for $x\leq z$. Thereafter, none of these markers can satisfy
the condition $s\leq d_s(x)\lor d_s(x)\neq d_{s-1}(x)$ for an
ordinary action. A coding action that affects $z$ must
have rank at most $z$, so $D$ also computes a stage beyond all such
coding actions. Taking a stage beyond these bounds and $z$ decides
whether $z\in A$. Conversely, if $n$ enters $D$ at stage $s$, then
after its coding action at most $n+1$ current markers remain in
$[0,s]$. The same bound holds if rank $n$ is absent. Subsequent
actions cannot increase this number, so $c_{n+1}>s$. Thus
$D(n)=D_{c_{n+1}}(n)$,
which proves $D\leq_TA$.
\end{proof}

We now find a threshold in $\BLR(A)$ above which all elements of
$\overline A$ have the same $e$-state. Fix $e$ and $q<2^{e+1}$.
We record the numbers that appear as markers with $e$-state at least
$q$. More precisely, at each stage $s$, before any coding or ordinary
action, the current markers are exactly
$[0,s]\setminus A_{s-1}$. Enumerate
each of these markers $y$ with $\nu_e(y,s)\geq q$ into $R_{e,q}$.
To use the domination property \eqref{eq:supermaximal-domination},
we also need a recursive search for a stage at which a suitable
marker is present. On input $x$, let $\Phi_{a(e,q)}$ simulate the
construction until the first stage $s\geq x$ at which the initial
marker list $[0,s]\setminus A_{s-1}$ contains a number $y>x$ with
$\nu_e(y,s)\geq q$. The function outputs this stage $s$; if no such
stage exists, it is undefined. The following \Cref{claim3} is immediate.

\begin{claim}\label{claim3}
The sets $R_{e,q}$ are uniformly r.e., and the indices $a(e,q)$ can
be chosen recursively. If $R_{e,q}$ is infinite, then
$\Phi_{a(e,q)}$ is total.
\end{claim}
Define $\nu_e(x)=\lim_s\nu_e(x,s)$.
\begin{claim}\label{claim4}
For every $e$ and $q<2^{e+1}$, if $R_{e,q}$ is infinite, $x\in\overline A$, $x>c_e$ and $x\geq h(a(e,q))$, then $\nu_e(x)\geq q$.
\end{claim}

\begin{proof}[Proof of the claim]
Suppose otherwise, and put
$s_0=\Phi_{a(e,q)}(x)$. Then
$s_0\leq d(x)$ by \eqref{eq:supermaximal-domination}. At stage
$s_0$ there is a marker $y>x$ of $e$-state at least $q$, whereas
$\nu_e(x,s)<q$ at every stage. Since $x>c_e$ and its rank cannot
increase, its current rank is greater than $e$ at every stage at
which it is present.

As long as $x$ survives, a marker greater than $x$ of $e$-state at
least $q$ continues to exist. Indeed, a tail action beginning
after $x$ either leaves such a witness untouched or replaces it by a
marker whose $e$-state is at least as large. An action beginning
at or before $x$ cannot retain $x$: if its rank is at least $e$,
maximality of the retained state rules this out; if its rank is less
than $e$, retaining $x$ would give it a rank less than $e$,
contrary to $x>c_e$.

If $d_{s_0}(x)\geq s_0$, then $x$ satisfies the condition
$s\leq d_s(x)\lor d_s(x)\neq d_{s-1}(x)$ at stage $s_0$. Otherwise,
$d_{s_0}(x)<s_0\leq d(x)$, so $d_s(x)$ changes at a later stage,
when the same condition is satisfied. At that stage, the persistent
witness makes $x$ eligible for an ordinary action. The action
at the least eligible rank therefore begins at or before $x$ and
cannot retain it, a contradiction. This proves
the claim.
\end{proof}

\begin{claim}
There is a function $b\in\BLR(D)$ such that, for every $e$ and
$q<2^{e+1}$, if $R_{e,q}$ is finite, then $R_{e,q}\subseteq[0,b(e,q)]$.
\end{claim}

\begin{proof}[Proof of the claim]
For an approximation with at most one change, initially output $0$
and use $\emptyset'$ to search for a number $b$ satisfying
$R_{e,q}\subseteq[0,b]$. On finding one, output it forever. If
there is no such bound, the approximation stays at $0$.
This defines $b\in\BLR(\emptyset')=\BLR(D)$ with the required property.
\end{proof}

\begin{claim}
\label{clm:supermaximal-threshold}
The function
\[
  f(e)=\max\bigl(\{c_e\}\cup
  \{h(a(e,q)),b(e,q):q<2^{e+1}\}\bigr)
\]
belongs to $\BLR(A)$ and witnesses that $A$ is supermaximal.
\end{claim}

\begin{proof}[Proof of the claim]
Since the function $e\mapsto c_e$ is $D$-recursive and the maximum
has a recursively bounded number of terms, their bounded-change approximations give
$f\in\BLR(D)=\BLR(A)$.
Fix $e$, $q<2^{e+1}$, and $x\in\overline A$ with $x>f(e)$.
Suppose first that $R_{e,q}$ is infinite. Since $x>c_e$ and
$x\geq h(a(e,q))$, \Cref{claim4} gives $\nu_e(x)\geq q$.
Now suppose that $R_{e,q}$ is finite. If $\nu_e(x)\geq q$, then
the permanent marker $x$ is eventually recorded in $R_{e,q}$.
But $R_{e,q}\subseteq[0,b(e,q)]$ and $x>b(e,q)$, a contradiction.
Hence $\nu_e(x)<q$ in this case.
Whether $R_{e,q}$ is infinite is independent of $x$, so all
elements of $\overline A$ greater than $f(e)$ have the same
$e$-state. In particular, they agree on membership in $W_e$.
Thus $A$ is supermaximal.
\end{proof}

Claims~\ref{clm:supermaximal-degree} and
\ref{clm:supermaximal-threshold} complete the proof.
\end{proof}

As a corollary, supermaximality is strictly stronger than maximality since there are high r.e. degrees that
are not superhigh~\cite{Mohrherr86}.

The density version of Martin's theorem,
\Cref{thm:density-martin}, shows that coarse domination characterizes highness.
The following proposition characterizes superhighness by requiring effective
control of how quickly the density of successful domination approaches one,
or of when it stays above any fixed positive density.
For a total function $f$ and an index $e$, put
\begin{equation*}
  D_e^f=\{n:\Phi_e(n)\mathord\downarrow
  \text{ and }\Phi_e(n)\leq f(n)\}.
\end{equation*}

\begin{proposition}
\label{thm:superhigh-density}
Fix a rational number $\alpha$ with $0<\alpha<1$. For every set $A$,
the following are equivalent.
\begin{enumerate}[label=\textnormal{(\roman*)}]
\item $A$ is superhigh.
\item There are a total function $f\leq_TA$ and a function
$m\colon\N^2\to\N$ with $m\in\BLR(A)$ such that, for every
$e\in\Tot$, every $k\geq1$, and every $N\geq\max\{1,m(e,k)\}$,
$\rho_N(D_e^f)\geq1-2^{-k}$.
\item There are a total function $f\leq_TA$ and a function
$h\in\BLR(A)$ such that, for every $e\in\Tot$ and every
$N\geq\max\{1,h(e)\}$,
$ \rho_N(D_e^f)\geq\alpha$.
\end{enumerate}
\end{proposition}

\begin{proof}
\textnormal{(i)}$\Rightarrow$\textnormal{(ii)}.
By \Cref{cor:clean-threshold}, there are
$f\leq_TA$ and $h\in\BLR(A)$ such that, for every
$e\in\Tot$, $(\forall n\geq h(e))\;\Phi_e(n)\leq f(n)$.
Define $m(e,k)=2^k(h(e)+1)$.
Then $m\in\BLR(A)$.
Fix $e\in\Tot$, $k\geq1$, and $N\geq\max\{1,m(e,k)\}$. Every member of
$[0,N)\setminus D_e^f$ is less than $h(e)$, so
$1-\rho_N(D_e^f)
  \leq\frac{h(e)}{N}
  \leq\frac{h(e)}{2^k(h(e)+1)}
  <2^{-k}$.

\textnormal{(ii)}$\Rightarrow$\textnormal{(iii)} is immediate.

\textnormal{(iii)}$\Rightarrow$\textnormal{(i)}.
Suppose that $f\leq_TA$ and $h\in\BLR(A)$ satisfy
\textnormal{(iii)}. Choose an integer $L\geq2$ with
$1/L<\alpha$. Fix a recursive function $\theta$ such that, whenever
$\Phi_e$ is total, $\Phi_{\theta(e)}(u)=\max_{i\leq u}\Phi_e(i)$
for every $u$. Define
\[
  F(n)=
  \begin{cases}
    \max\{f(u):n\leq u<Ln\},&\text{if }n\geq1,\\
    f(0),&\text{if }n=0,
  \end{cases}
\]
and
$h_0(e)=\max\left\{1,
  \left\lceil\frac{h(\theta(e))}{L}\right\rceil\right\}$.
Then $F\leq_TA$ and $h_0\in\BLR(A)$. 
If $e\in\Tot$ and $n\geq h_0(e)$, then $Ln\geq h(\theta(e))$.
The interval $[n,Ln)$ must meet $D_{\theta(e)}^f$, since otherwise
$\rho_{Ln}(D_{\theta(e)}^f)\leq\frac1L<\alpha$.
For $u$ in this intersection,
\[
  F(n)\geq f(u)\geq\Phi_{\theta(e)}(u)
  =\max_{i\leq u}\Phi_e(i)\geq\Phi_e(n).
\]
Thus $F$ and $h_0$ satisfy \Cref{cor:clean-threshold}, and $A$ is
superhigh.
\end{proof}

\section{The non-robustness of pdomination under lower density}
\label{sec:ldpd-collapse}

Pdomination was introduced by Chong, Hoi, Stephan, and
Turetsky~\cite{ChongHoiStephanTuretsky15}. The notion was
motivated by the bounded monotone enumeration principle BME in reverse mathematics; see \cite{CSY}.

\begin{definition}[\cite{ChongHoiStephanTuretsky15}]
\label{def:pdominance}
A set $A$ is \emph{pdominant} if there is a partial $A$-recursive function
$\psi$ such that, for every partial recursive $\varphi$,
\[
  (\forall^\infty x\in\operatorname{dom}(\varphi))
  (\exists y\leq x)\,[\psi(y)\mathord\downarrow>\varphi(x)].
\]
\end{definition}

Chong, Hoi, Stephan, and Turetsky showed that pdomination is a nontrivial
degree property.  No weakly $2$-generic set is pdominant, whereas there are pdominant
$1$-generic sets below the halting problem; the two halves of
Chaitin's $\Omega$ are pdominant, but no set which is low for Martin-L\"of
randomness is pdominant; there are a low r.e. pdominant set and a
high r.e. set which is not pdominant.

We consider two density weakenings of pdomination, using correctness sets
and density relative to the test domain, respectively.  We show that both
versions hold for every set, first at positive lower density and then at
density one.

For partial functions $\psi,\varphi$,
put
\begin{equation*}
  H(\psi,\varphi)=
  \{x\in\operatorname{dom}(\varphi):
  (\exists y\leq x)\,
  [\psi(y)\mathord\downarrow>\varphi(x)]\},
\end{equation*}
and
\begin{equation*}
  C(\psi,\varphi)=
  (\N\setminus\operatorname{dom}(\varphi))\cup H(\psi,\varphi).
\end{equation*}

\begin{definition}\label{def:correctness-pdominance}
A set $A$ is \emph{lower-correctness-pdominant}, abbreviated \emph{lcPD}, if
there is a partial $A$-recursive function $\psi$ such that, for every partial
recursive $\varphi$ with infinite domain,
$\underline\rho
  (C(\psi,\varphi))>0$.
Replacing $>0$ by $=1$ gives \emph{correctness-pdominance}, abbreviated
\emph{cPD}.
\end{definition}

For an
infinite set $D\subseteq\N$ and $S\subseteq\N$, put
$\underline\rho_D(S)=\liminf_{N\to\infty}
  \frac{|S\cap D\cap[0,N)|}{|D\cap[0,N)|}$.
When the corresponding limit exists, we write
\[
  \rho_D(S)=\lim_{N\to\infty}
  \frac{|S\cap D\cap[0,N)|}{|D\cap[0,N)|}.
\]

\begin{definition}\label{def:lower-density-pdominance}
A set $A$ is \emph{lower-density-pdominant}, abbreviated \emph{ldPD}, if there
is a partial $A$-recursive function $\psi$ such that, for every partial
recursive $\varphi$ with infinite domain,
$\underline\rho_{\operatorname{dom}(\varphi)}
  (H(\psi,\varphi))>0$.
Replacing $>0$ by $=1$ gives \emph{density-pdominance}, abbreviated
\emph{dPD}.
\end{definition}
Every pdomination witness is a dPD witness.  Every dPD witness is both an
ldPD witness and a cPD witness, and every cPD witness is an lcPD witness.

We first give a single partial recursive witness for lcPD and ldPD.

\begin{theorem}
\label{prop:uniform-positive-density-degeneration}
There is a single partial recursive function $\psi$ such that, for
every $e\in\N$,
$\underline\rho\bigl(C(\psi,\Phi_e)\bigr)
  \geq \frac{1}{2^{e+3}}>0$.
Moreover, if $D=\operatorname{dom}(\Phi_e)$ is infinite, then
$\underline\rho_D\bigl(H(\psi,\Phi_e)\bigr)
  \geq \frac{1}{2^{e+3}}>0$.
Thus every set is both \emph{lcPD} and \emph{ldPD}.
\end{theorem}

\begin{proof}
The idea is to partition the inputs for each test
$\Phi_e$ into cells of a fixed length and to force one genuine hit whenever
a cell contains a convergence.  A dedicated marker is placed below every
point of its cell; once the first convergence in that cell is seen, we define
$\psi$ at the marker late enough to exceed the observed value of $\Phi_e$. Disjoint marker columns allow this to
be done for all $e$ simultaneously.
Thus every nonempty cell contributes one hit but contains at most a fixed
number of convergence points, which gives the required positive-density
bound.

For the requirement indexed by $e$, put
$L=L_e=2^{e+3}$,
and divide a tail of $\N$ into cells
$ I_{e,j}=[jL,(j+1)L)$, $j\geq1$.
Assign to the $j$th cell the marker
$a_{e,j}=2^{e+1}(2j+1)-1$.
The marker columns belonging to distinct $e$ are disjoint, since
$2^{e+1}$ is the exact power of $2$ dividing $a_{e,j}+1$.
Moreover, if $x\in I_{e,j}$, then
$a_{e,j}<j2^{e+3}=jL\leq x$.

Construct an r.e. set $W$ as follows: for all
$e,j$, when a convergence $\Phi_e(x)\mathord\downarrow=t$ with
$x\in I_{e,j}$ is first found, wait until a construction stage $s>t$ and
enumerate $a_{e,j}$ into $W$.

Each marker is enumerated at most once.  Define
$\psi(a)=\min\{s:a\in W_s\}$
when $a\in W$, and leave $\psi(a)$ undefined otherwise.  This is a partial
recursive function.  Designate one convergence $x_{e,j}$ witnessing the
action of each cell.  For every designated $x=x_{e,j}$, $a_{e,j}\leq x$ and $ \psi(a_{e,j})>\Phi_e(x)$, so
$x\in H(\psi,\Phi_e)$.

Then it is clear that $\underline\rho\bigl(C(\psi,\Phi_e)\bigr) \ge \underline\rho_D\bigl(H(\psi,\Phi_e)\bigr)
  \geq \frac{1}{2^{e+3}}$ if $D$ is infinite. And if $D$ is finite, $C(\psi,\Phi_e)$ is cofinite and has density 1.
\end{proof}

To obtain density one, we use several markers in each cell and repeat the
construction at every prescribed accuracy.  The following lemma supplies
the thresholds for these markers.

\begin{lemma}
\label{lem:finite-multiplicative-net}
Given a rational number $0<\varepsilon<1$ and an integer $M\geq1$, one can
effectively find a multiple $L$ of $M$ and a finite set
$T\subseteq\{1,\ldots,L\}$ such that
$|T|\leq\frac{L}{M}$
and, for every integer $1\leq r\leq L$, there is some $t\in T$ satisfying
$ (1-\varepsilon)r\leq t\leq r$.
\end{lemma}

\begin{proof}
For any $L\geq1$, define $t_0=1$,
$t_{i+1}=\left\lceil\frac{t_i}{1-\varepsilon}\right\rceil$,
let $m$ be least such that $t_m>L$, and put
$T_L=\{t_0,\ldots,t_{m-1}\}$.  Given $1\leq r\leq L$, choose the largest
$i<m$ such that $t_i\leq r$.  Then $r<t_{i+1}$, and, since $r$ is an integer,
$r\leq t_{i+1}-1
  <\frac{t_i}{1-\varepsilon}$.
Consequently $(1-\varepsilon)r<t_i\leq r$, so $T_L$ has the required
covering property.

Moreover, $t_i\geq(1-\varepsilon)^{-i}$ for every $i$.  Since
$t_{m-1}\leq L$,
\[
  |T_L|=m
  \leq 1+\frac{\log L}{-\log(1-\varepsilon)}
  \leq 1+\varepsilon^{-1}\log L,
\]
where the last inequality follows from
$-\log(1-\varepsilon)\geq\varepsilon$.  Hence, for all sufficiently large
multiples $L$ of $M$,
$|T_L|\leq\frac{L}{M}$.

To find such an $L$ effectively, search through the positive multiples of
$M$, construct the corresponding finite set $T_L$, and stop when $|T_L|\leq L/M$.
The preceding estimate guarantees that this search terminates.  Taking $T=T_L$ completes the proof.
\end{proof}

\begin{theorem}
\label{thm:dpd-collapse}
Every set is both \emph{cPD} and \emph{dPD}.  Indeed, there is a  partial
recursive function $\psi$ such that, for every partial recursive function
$\varphi$ with infinite domain, putting
$D=\operatorname{dom}(\varphi)$, $\rho\bigl(C(\psi,\varphi)\bigr)=1$ and $\rho_D\bigl(H(\psi,\varphi)\bigr)=1$.
\end{theorem}

Since the final number $r$ of convergence points in a cell is not known while
the computation unfolds, at each prescribed accuracy $\varepsilon_n$ we
reserve several markers indexed by the thresholds
from \Cref{lem:finite-multiplicative-net}.  When the
$t$th convergence appears, the corresponding marker makes the first $t$
convergence points into hits; once all convergences in the cell have appeared,
some threshold $t$
lies between $(1-\varepsilon_n)r$ and $r$.  Hence all but an
$\varepsilon_n$-fraction of the convergence points in each complete cell are
hits, while disjoint reservoirs let all pairs $(e,n)$ act simultaneously.

\begin{proof}
Fix $e,n$ and put $\varepsilon_n=2^{-n-1}$ and $q=\langle e,n\rangle$.
Use the disjoint recursive marker reservoir $R_q=\{2^q(2r+1)-1:r\in\N\}$.
Apply \Cref{lem:finite-multiplicative-net} with
$\varepsilon=\varepsilon_n$ and $M=2^{q+1}$ to choose effectively a multiple $L=L_{e,n}$ of $2^{q+1}$
and a finite set
$T\subseteq\{1,\ldots,L\}$ such that $|T|\leq L/2^{q+1}$ and, for every
$1\leq r\leq L$, some $t\in T$ satisfies $(1-\varepsilon_n)r\leq t\leq r$.

For $j\geq0$, let
$I_{q,j}=[(j+1)L,(j+2)L)$.
The preceding interval $[jL,(j+1)L)$ contains exactly $L/2^{q+1}$ members of
$R_q$, so assign distinct such markers $a_{q,j,t}$ to the thresholds $t\in
T$.  When $t$ distinct convergences of $\Phi_e$ in $I_{q,j}$ have been
seen, define
\[
  \psi(a_{q,j,t})=1+\max\{\Phi_e(x):x\text{ is one of these }t
  \text{ convergence points}\}.
\]
Dovetailing over all $e,n,j,t$ defines a partial recursive function, since
the reservoirs $R_q$ are pairwise disjoint.  Moreover,
$a_{q,j,t}<\min I_{q,j}$, so the first $t$ convergence points in the cell
belong to $H(\psi,\Phi_e)$.

Fix $e,n$, put $D=\operatorname{dom}(\Phi_e)$, and let
$E=D\setminus H(\psi,\Phi_e)$.  We first spell out the estimate on a
complete cell.  Put
$r_j=|D\cap I_{q,j}|$.
If $r_j>0$, choose $t_j\in T$ so that
$(1-\varepsilon_n)r_j\leq t_j\leq r_j$.
The inequality $t_j\leq r_j$ ensures that the construction eventually sees
$t_j$ distinct convergence points of $\Phi_e$ in $I_{q,j}$ and defines
$\psi(a_{q,j,t_j})$.  Each of these $t_j$ points lies in
$H(\psi,\Phi_e)$.  Thus $|H(\psi,\Phi_e)\cap I_{q,j}|\geq t_j$
and consequently
\begin{equation}\label{eq:dpd-cell-error}
  |E\cap I_{q,j}|
  =r_j-|H(\psi,\Phi_e)\cap I_{q,j}|
  \leq r_j-t_j
  \leq\varepsilon_n r_j.
\end{equation}
When $r_j=0$, the same final estimate is immediate.

Now fix $N$.  The cells $I_{q,j}$ ($j \in \omega$) are consecutive and cover $[L,\infty)$.
Hence $[0,N)$ consists of the initial piece
$[0,\min\{L,N\})$, all the cells $I_{q,j}$ that are wholly contained in
$[0,N)$, and possibly an initial segment of one further cell.  The initial piece and the terminal incomplete cell each
have length at most $L$, so together they contain at most $2L$ points of
$E$.  Applying \eqref{eq:dpd-cell-error} to every complete cell and summing
gives
\begin{align*}
  |E\cap[0,N)|
  &\leq 2L+
    \sum_{I_{q,j}\subseteq[0,N)}|E\cap I_{q,j}|\\
  &\leq 2L+\varepsilon_n
    \sum_{I_{q,j}\subseteq[0,N)}|D\cap I_{q,j}|\\
  &\leq 2L+\varepsilon_n|D\cap[0,N)|.
\end{align*}
If $D$ is infinite, then $|D\cap[0,N)|\to\infty$. Then $\limsup_{N\to\infty}
  \frac{|E\cap[0,N)|}{|D\cap[0,N)|}\leq\varepsilon_n$.
Since this holds for every $n$, the limsup is zero, and $\rho_D\bigl(H(\psi,\Phi_e)\bigr)=1$.  Finally, since $\underline\rho\bigl(C(\psi,\Phi_e)\bigr) \ge \underline\rho_D\bigl(H(\psi,\Phi_e)\bigr)$, we have $\rho\bigl(C(\psi,\Phi_e)\bigr)=1$ as well.
\end{proof}

\section*{acknowledgments}
The author was supported by Focused Research Grants of Natural Science Foundation of Jiangsu Province No. BK20243060. The author thanks Professor Liang Yu for helpful discussions, and Theorem \ref{thm:superhigh-threshold} was joint work with him.
The assertion of Theorem \ref{thm:superhigh-supermaximal} was posed as a question by the author in a student talk given at NUS when he was supported by the Institute for Mathematical Sciences to attend the 2026 logic summer school. The proof of this Theorem was later produced by ChatGPT 6. This proof has been verified by the author without the use
of AI tools. The author takes full responsibility for the final content of the paper.

\end{document}